\documentclass{article}[12pt]
\usepackage{graphicx,pstricks,comment}
\usepackage[all]{xy}
\usepackage{subfigure}
\usepackage{rotating}
\usepackage{amsfonts}
\usepackage{amssymb}
\usepackage{amsmath}
\usepackage{xcolor}
\usepackage[utf8]{inputenc}
\usepackage{tikz}
\usepackage{pbox}

\usepackage{indentfirst}

\usepackage{epsfig}
\usepackage{psfrag}
\usepackage{enumerate}
\usepackage{array}
\usepackage{theorem}

\newcommand{\C}{\mathbb{C}}
\newcommand{\R}{\mathbb{R}}

\newcommand{\Z}{\mathbb{Z}}

\newcommand{\Q}{\mathbb{Q}}

\newcommand{\HCn}{{\rm H}_\C^n}

\newcommand{\PU}{{\rm PU}}
\newcommand{\SU}{{\rm SU}}

\newtheorem{thm}{Theorem}
\newtheorem{cor}{Corollary}
\newtheorem{lem}{Lemma}
\newtheorem{prop}{Proposition}

{\theorembodyfont{\rmfamily} }
{\theorembodyfont{\rmfamily} }

\newcommand{\Pf}{{\em Proof}. }
\newcommand{\EPf}{\hfill$\Box$\vspace{.5cm}}

\title{Hybrid lattices and thin subgroups of Picard modular groups}

\author{Julien Paupert\footnote{Both authors partially supported by National Science Foundation Grant DMS-1708463.}, Joseph Wells$^*$}

\begin{document}
\maketitle

\begin{abstract} We consider a certain hybridization construction which produces a subgroup of $\PU(n,1)$ from a pair of lattices in $\PU(n-1,1)$. Among the Picard modular groups $\PU(2,1,\mathcal{O}_d)$, we show that the hybrid of pairs of Fuchsian subgroups $\PU(1,1,\mathcal{O}_d)$ is a lattice when $d=1$ and $d=7$, and a geometrically infinite thin subgroup when $d=3$, that is an infinite-index subgroup with the same Zariski-closure as the full lattice.  

\end{abstract}


\section{Introduction}

Lattices in rank 1 real (semi)simple Lie groups are still far from understood. A key notion is that of \emph{arithmetic lattice} which we will not define properly here but note that by a famous result of Margulis a lattice in such a Lie group is arithmetic if and only if it has infinite index in its commensurator.

Margulis' celebrated arithmeticity theorem states that every lattice of a simple real Lie group $G$ is arithmetic whenever the real rank of $G$ is at least two. Thus non-arithmetic lattices can only exist in real rank one, that is when the associated symmetric space is a hyperbolic space. In real hyperbolic space, where the Lie group is ${\rm PO}(n,1)$, Gromov and Piatetski-Shapiro  produced in \cite{GPS} a construction yielding non-arithmetic lattices in ${\rm PO}(n,1)$ for all $n \geqslant 2$ (see below for more details), in fact producing in each dimension infinitely many non-commensurable lattices, both cocompact and non-cocompact.
In quaternionic hyperbolic spaces (and the Cayley hyperbolic plane), work of Corlette and Gromov-Schoen implies as in the higher rank case that all lattices are arithmetic.

The case of complex hyperbolic spaces, where the associated Lie group is $\PU(n,1)$, is much less understood. Non-arithmetic lattices in $\PU(2,1)$ were first constructed by Mostow in 1980 in \cite{M1}, and subsequently by Deligne-Mostow and Mostow as monodromy groups of certain hypergeometric functions in \cite{DM} and \cite{M2}, following pioneering work of Picard. More recently, Deraux, Parker and the first author constructed new families of non-arithmetic lattices in $\PU(2,1)$ by considering groups generated by certain triples of complex reflections (see \cite{DPP1}, \cite{DPP2}).  
Taken together, these constructions yield 22 commensurability classes of non-arithmetic lattices in $\PU(2,1)$, and only 2 commensurability classes in $\PU(3,1)$. The latter two are noncocompact; one is a Deligne-Mostow lattice and the other was constructed by Couwenberg-Heckman-Looijenga in 2005 and recently found to be non-arithmetic by Deraux, \cite{D}.
Major open questions in this area remain the existence of non-arithmetic lattices in $\PU(n,1)$ for $n \geqslant 4$, as well as the number (or finiteness thereof) of commensurability classes in each dimension. 

The Gromov--Piatetski-Shapiro construction, which they call interbreeding of 2 arithmetic lattices (now often referred to as hybridization), produces a lattice $\Gamma < {\rm PO}(n,1)$ from 2 lattices $\Gamma_1$ and $\Gamma_2$ in  ${\rm PO}(n,1)$ which have a common sublattice $\Gamma_{12} < {\rm PO}(n-1,1)$. Geometrically, this provides two hyperbolic $n$-manifolds $V_1=\Gamma_1 \backslash H_\R^n$ and $V_2=\Gamma_2 \backslash H_\R^n$ with a hyperbolic $(n-1)$-manifold $V_{12}$ which is isometrically embedded in $V_1$ and $V_2$ as a totally geodesic hypersurface. This allows one to produce the hybrid manifold $V$ by gluing $V_1-V_{12}$ and $V_2-V_{12}$ along $V_{12}$ (more precisely, in case $V_{12}$ separates $V_1$ and $V_2$, by gluing $V_1^+ -V_{12}$ and $V_2^+ -V_{12}$ along $V_{12}$, with $V_i^+$ a connected component of $V_i - V_{12}$). The resulting manifold is also hyperbolic because the gluing took place along a totally geodesic hypersurface, and its fundamental group $\Gamma$ is therefore a lattice in ${\rm PO}(n,1)$. The main point is then that if $\Gamma_1$ and $\Gamma_2$ are both arithmetic but non-commensurable, their hybrid $\Gamma$ is non-arithmetic. Note that the resulting hybrid $\Gamma$ is algebraically an amalgamated free product of $\Gamma_1$ and $\Gamma_2$ over $\Gamma_{12}$ (say, in the case where $V_{12}$ separates both $V_1$ and $V_2$), and in all cases is generated by its sublattices $\Gamma_1$ and $\Gamma_2$.

It is not straightforward to adapt this construction to construct lattices in $\PU(n,1)$, the main difficulty being that there do not exist in complex hyperbolic space any totally geodesic real hypersurfaces. In fact, it has been a famous open question since the work of Gromov--Piatetski-Shapiro to find some analogous construction in $\PU(n,1)$. Hunt proposed the following construction (see \cite{Pau} and references therein). Start with 2 arithmetic lattices $\Gamma_1$ and $\Gamma_2$ in $\PU(n,1)$, and suppose that one can embed them in $\PU(n+1,1)$ in such a way that (a) each stabilizes a totally geodesic $H_\C^n \subset H_\C^{n+1}$ (b) these 2 complex hypersurfaces are orthogonal, and (c) the intersection of the embedded $\Gamma_i$ is a lattice in the corresponding $\PU(n-1,1)$. The resulting hybrid $\Gamma=H(\Gamma_1,\Gamma_2)$ is then defined as the subgroup of $\PU(n+1,1)$ generated by the images of $\Gamma_1$ and $\Gamma_2$. (See the end of Section~\ref{chg} for a more detailed and concrete description when $n=2$).

It is not clear when, if ever, such a group has any nice properties. One expects in general the hybrid group to be non-discrete, and in fact the first author showed in \cite{Pau} that this happens infinitely often among hybrids in $\PU(2,1)$ of pairs of Fuchsian triangle subgroups of $\PU(1,1)$. It was observed there that one can easily arrange for the hybrid to be discrete by arranging for the two subgroups $\Gamma_1,\Gamma_2$ to already belong to a known lattice. But even in the simplest case of arithmetic cusped lattices (where the matrix entries are all in $\mathcal{O}_d$, the ring of integers of $\Q[i\sqrt{d}]$ for some squarefree $d \geqslant 1$), it was not known whether the discrete hybrid $\Gamma$ could ever be a sublattice of the corresponding Picard modular group $\Gamma(d)=\PU(2,1,\mathcal{O}_d)$, as opposed to an infinite-index (discrete) subgroup of $\Gamma(d)$. Following Sarnak (\cite{S}) we will call \emph{thin subgroup} of a lattice $\Gamma$ any infinite-index subgroup having the same Zariski-closure as $\Gamma$.

In this note we show that in fact both behaviors can occur, even among this simplest class of hybrids of sublattices of the Picard modular groups $\Gamma(d)$. More precisely, we consider for $d=3,1,7$ the hybrid subgroup $H(d)$ defined as the hybrid of two copies of $\SU(1,1,\mathcal{O}_d)$ inside the Picard modular group $\PU(2,1,\mathcal{O}_d)$. 
These specific values of $d$ are those for which a presentation of $\PU(2,1,\mathcal{O}_d)$ is known (by \cite{FP}, \cite{FFP} and \cite{MP}). Our main results can be summarized as follows (combining Theorems~\ref{d=3main}, \ref{d=1main}, and \ref{d=7main} and Propositions~\ref{limitsets} and \ref{geominf}).

\begin{thm}\label{main} 
\begin{enumerate}
\item The hybrid $H(3)$ is a thin subgroup of the Eisenstein-Picard lattice $\PU(2,1,\mathcal{O}_3)$. It has full limit set $\partial_\infty {\rm H^2_\C} \simeq S^3$ and is therefore geometrically infinite.
\item The hybrid $H(1)$ has index 2 in the Gauss-Picard lattice $\PU(2,1,\mathcal{O}_1)$.
\item The hybrid $H(7)$ is the full Picard lattice $\PU(2,1,\mathcal{O}_7)$.
\end{enumerate}
\end{thm}

{\bf Remarks:} 

\smallskip

(a) We also give analogous results for two related hybrids $H'(3)$ and $H'(1)$ in Corollaries~\ref{d=3other} and \ref{d=1other}. In terms of Fuchsian triangle groups these groups are defined as the hybrids of two copies of the (orientation-preserving) triangle groups $(2,6,\infty)$ and $(2,4,\infty)$ respectively, as opposed to  $(3,\infty,\infty) \simeq \PU(1,1,\mathcal{O}_3)$ and $(2,\infty,\infty) \simeq \PU(1,1,\mathcal{O}_1)$ (so, replacing the elliptic generator by one of its square roots).  
An interesting feature of $H'(3)$ is that it has infinite index in its normal closure in $\Gamma(3)$, whereas all other hybrids we consider are normal in $\Gamma(d)$.

(b) In all cases we also show that the hybrid $\Gamma$ is not an amalgamated free product of $\Gamma_1$ and $\Gamma_2$ over their intersection. In case $\Gamma$ is itself a lattice this follows from general considerations of cohomological dimension, and for $H(3)$ and $H'(3)$ we show this by finding sufficiently many relations among the generators for $\Gamma$, see Corollary~\ref{notamalgam}.

(c) One of the main geometric difficulties in analyzing these groups is understanding the parabolic subgroups. By construction the generators contain a pair of (opposite) parabolic isometries (as well as an elliptic isometry when $d=3$, two elliptic isometries when $d=1$, and two elliptic and two loxodromic isometries when $d=7$), however it seems hard in general to determine the rank of the parabolic subgroups of the hybrid. In the cases where the hybrid is a lattice we obtain indirectly that the parabolic subgroups must have full rank, but in the thin subgroup case we do not know what this rank is.  

(d) The parabolic isometries appearing in the generators for our hybrids are by construction vertical Heisenberg translations, since they preserve a complex line (see Section~\ref{chg}). It turns out that Falbel (\cite{F}) and Falbel-Wang (\cite{FW}) studied a group formally similar to our hybrid $H(3)$, obtained by completely different methods, namely by finding all irreducible representations of the figure-eight knot group $\Gamma_8$ into $\PU(2,1)$ with unipotent boundary holonomy. Falbel showed in \cite{F} that there are exactly 3 such representations, one of which has image contained in  $\Gamma(3)=\PU(2,1,\mathcal{O}_3)$ and the two others in  $\Gamma(7)=\PU(2,1,\mathcal{O}_7)$. These are all generated by a pair of opposite horizontal Heisenberg translations. The image of the former representation is shown in \cite{F} and \cite{FW} to be, like our hybrids $H(3)$ and $H'(3)$, a thin subgroup of $\Gamma(3)$ with full limit set, whereas the images of the latter two representations are shown in \cite{DF} to have non-empty domain of discontinuity (and hence have infinite index in $\Gamma(7)$). We were inspired by some of the arguments of \cite{F} and \cite{FW}.

(e) Kapovich found in \cite{K1} the first examples of infinite-index normal subgroups of lattices in $\PU(2,1)$, among a family of four lattices first constructed by Livn\'e in his thesis (and predating the term \emph{thin subgroup}). Parker showed in \cite{Par0} (sections 6 and 7) that this description could be extended to the Eisenstein-Picard modular group $\PU(2,1,\mathcal{O}_3)$, and that  some of Kapovich's results extended to that case as well. An anonymous referee showed us by a sequence of clever computations that our hybrid $H(3)$ is in fact commensurable to the infinite-index normal subgroup of $\PU(2,1,\mathcal{O}_3)$ obtained in this way. It would be interesting to understand this relation geometrically.

(f) Discrete groups generated by opposite parabolic subgroups have been studied in higher rank by Oh, Benoist-Oh and others. A conjecture of Margulis states that if $G$ is a semisimple real algebraic group of rank at least 2 and $\Gamma$ a discrete Zariski-dense subgroup containing irreducible lattices in two opposite horospherical subgroups, then $\Gamma$ is an arithmetic lattice in $G$. Oh showed in \cite{O} that this holds when $G$ is a simple split real Lie group except ${\rm SL}(3,\R)$,  Benoist-Oh extended this in \cite{BO} to the case of $G={\rm SL}(3,\R)$, and very recently Benoist-Miquel treated the general case in \cite{BM}.

\smallskip

The paper is organized as follows. In section 2 we review basic facts about complex hyperbolic space, its isometries, subspaces and boundary at infinity. In Sections 3,4,5 we consider each of the hybrids $H(3)$, $H(1)$ and $H(7)$ respectively. In section 6 we review and apply basic facts about limit sets and geometrical finiteness to the non-lattice hybrid $H(3)$. We would like to thank Elisha Falbel for pointing out a simplification of the proof of Theorem~\ref{d=3main}, and an anonymous referee for several useful comments.

\section{Complex hyperbolic space, isometries and boundary at infinity}\label{chg}

We give a brief summary of basic definitions and facts about complex hyperbolic geometry, and refer the reader to \cite{G}, \cite{CG} or \cite{Par2} for more details.

\smallskip

{\bf Projective models of ${\rm H}_\C^n$:} 

\smallskip

Denote $\C^{n,1}$ the vector space $\C^{n+1}$ endowed with a Hermitian form $\langle \cdot \, , \cdot \rangle$ of signature $(n,1)$. Define $V^-=\left\lbrace Z \in \C^{n,1} | \langle Z , Z \rangle <0 \right\rbrace$ and $V^0=\left\lbrace Z \in \C^{n,1} | \langle Z , Z \rangle =0 \right\rbrace$.
Let $\pi: \C^{n+1}-\{0\} \longrightarrow \C{\rm P}^n$ denote projectivization.
One may then define complex hyperbolic $n$-space ${\rm H}_\C^n$ as $\pi(V^-) \subset \C{\rm P}^n$, with the distance $d$ (corresponding to the Bergman metric) given by:

\begin{equation}\label{dist}
  \cosh ^2 \frac{1}{2}d(\pi(X),\pi(Y)) = \frac{|\langle X, Y \rangle|^2}{\langle X, X \rangle  \langle Y, Y \rangle}
\end{equation}

The boundary at infinity $\partial {\rm H}_\C^n$ is then naturally identified with $\pi(V_0)$.
 Different Hermitian forms of signature $(n,1)$ give rise to different models of $\HCn$. Two of the most common choices are the Hermitian forms corresponding to the Hermitian matrices $H_1={\rm Diag}(1,...,1,-1)$ and:
  
\begin{equation}\label{H2matrix}
H_2=
\begin{bmatrix}
0 & 0 & 1 \\
0 & I_{n-1}& 0 \\
1 & 0 & 0 
\end{bmatrix}
\end{equation}
 
 In the first case, $\pi(V^-) \subset \C P^n$ is the unit ball of $\C^n$, seen in the affine chart $\{ z_{n+1}=1\}$ of $\C P^n$, hence the model is called the \textit{ball model} of $\HCn$. 
In the second case,  we obtain the \emph{Siegel model} of $\HCn$, which is analogous to the upper-half space model of ${\rm H}_\R^n$ and is likewise well-adapted to parabolic isometries fixing a specific boundary point. We will mostly use the Siegel model in this paper and will give a bit more details about it below. We will use the following \emph{Cayley transform} $J$ to pass from the ball model to the Siegel model (see \cite{Par2}); a key point for us is that $J \in {\rm GL}(3,\Z)$, hence conjugating by $J$ preserves integrality of matrix entries.
 
\begin{equation}\label{cayley}
J = \begin{pmatrix} 1 & 1 & 0 \\ 0 & 1 & -1 \\ 1 & 1 & -1 \end{pmatrix}
\end{equation}

{\bf Isometries:}

\smallskip

It is clear from \eqref{dist} that $\PU(n,1)$ acts by isometries on ${\rm H}_\C^n$, denoting ${\rm U}(n,1)$ the subgroup of ${\rm GL}(n+1,\C)$ preserving the Hermitian form, and $\PU(n,1)$ its image in ${\rm PGL}(n+1,\C)$. It turns out that PU($n$,1) is the group of holomorphic isometries of ${\rm H}_\C^n$, and the full group of isometries is $\PU(n,1) \ltimes \Z/2$, where the $\Z/2$ factor corresponds to a real reflection (see below). A holomorphic isometry of $\HCn$ is of one of the following three types:
\begin{itemize}
\item \emph{elliptic} if it has a fixed point in ${\rm H}_\C^n$
\item \emph{parabolic} if it has (no fixed point in ${\rm H}_\C^n$ and) exactly one fixed point in $\partial{\rm H}_\C^n$
\item \emph{loxodromic}: if it has (no fixed point in ${\rm H}_\C^n$ and) exactly two fixed points in $\partial{\rm H}_\C^n$
\end{itemize}

{\bf Totally geodesic subspaces:} 

\smallskip

 A \emph{complex k-plane} is a projective $k$-dimensional subspace of $\C P^n$ intersecting $\pi(V^-)$  non-trivially (so, it is an isometrically embedded copy of ${\rm H}_\C^{k} \subset {\rm H}_\C^n$). Complex 1-planes are usually called \emph{complex lines}. If $L=\pi(\tilde{L})$ is a complex $(n-1)$-plane, any $v \in \C^{n+1}-\{ 0\}$ orthogonal to $\tilde{L}$ is called a \emph{polar vector} for $L$. 
 
 A \emph{real k-plane} is the projective image of a totally real $(k+1)$-subspace $W$ of $\C^{n,1}$, i. e. a $(k+1)$-dimensional real linear subspace such that $\langle v,w \rangle \in \R$ for all $v,w \in W$.  We will usually call real 2-planes simply real planes, or $\R$-planes. Every real $n$-plane in ${\rm H}_\C^n$ is the fixed-point set of an antiholomorphic isometry of order 2 called a \emph{real reflection} or $\R$-reflection. The prototype of such an isometry is the map given in affine coordinates by $(z_1,...,z_n) \mapsto (\overline{z_1},...,\overline{z_n})$; this is an isometry provided that the Hermitian form has real coefficients.

We will need to distinguish between the following types of parabolic isometries. A parabolic isometry is called \emph{unipotent} if it has a unipotent lift to ${\rm U}(n,1)$. In dimensions $n>1$, unipotent isometries are either {\it 2-step} (also called \emph{vertical}) or {\it 3-step} (also called \emph{horizontal}), according to whether the minimal polynomial of their unipotent lift is $(X-1)^2$ or $(X-1)^3$ (see section 3.4 of \cite{CG}). Another way to distinguish these two types is that 2-step unipotent isometries preserve a complex line (in fact, any complex line through their fixed point) but no real plane, whereas 3-step unipotent isometries preserve a real plane (in fact, an entire \emph{fan} of these, see section 2.3 of \cite{PW}) but no complex line.  

\smallskip

{\bf Boundary at infinity and Heisenberg group:}

\smallskip

In the Siegel model associated to the Hermitian form given by the matrix $H_2$ in (\ref{H2matrix}), ${\rm H}_\C^n$ can be parametrized by $\C^{n-1} \times \R \times \R^+$ as follows, denoting as before by $\pi$ the projectivization map: ${\rm H}_\C^n=\{\pi(\psi(z,t,u)) \, | \, z \in \C^{n-1}, t \in \R, u \in \R^+\}$, where: 
\begin{eqnarray}\label{horocoord}
  \psi(z,t,u)= \left(\begin{array}{c}
    (-|z|^2-u+it)/2 \\
    z \\
    1
  \end{array}\right)
\end{eqnarray}
With this parametrization the boundary at infinity $\partial_\infty {\rm H}_\C^n$ corresponds to the one-point compactification:
$$\left\{\pi\!\left(\psi(z,t,0)\right) \, | \, z \in \C^{n-1}, t \in \R \right\} \cup \{\infty\}$$
where $\infty=\pi((1,0,...,0)^T)$. The coordinates $(z,t,u) \in \C^{n-1} \times \R \times \R^+$ are called the \emph{horospherical coordinates} of the point $\pi(\psi(z,t,u) \in {\rm H}_\C^n$. 

The punctured boundary  $\partial_\infty {\rm H}_\C^n - \{ \infty \}$ is then naturally identified to the \emph{generalized Heisenberg group} ${\rm Heis}(\C,n)$, defined as the set $\C^{n-1} \times \R$ equipped with the group law: $$(z_1,t_1)(z_2,t_2)=(z_1+z_2,t_1+t_2+2{\rm Im} \, (z_1 \cdot \overline{z_2}))$$ where $\cdot$ denotes the usual Euclidean dot-product on $\C^{n-1}$. This is the classical 3-dimensional Heisenberg group when $n=2$. The identification of  $\partial_\infty {\rm H}_\C^n - \{ \infty \}$ with ${\rm Heis}(\C,n)$ is given by the simply-transitive action of ${\rm Heis}(\C,n)$ on $\partial_\infty {\rm H}_\C^n - \{ \infty \}$, where the element $(z_1,t_1) \in {\rm Heis}(\C,n)$ acts on the vector $\psi(z_2,t_2,0)$ by left-multiplication by the following \emph{Heisenberg translation} matrix in ${\rm U}(n,1)$:
\begin{equation}\label{heistrans}
  T_{(z_1,t_1)}=\left(\begin{array}{ccc}
    1 & -z_1^* & (-|z_1|^2+it_1)/2 \\
    0 & {\rm I}_{n-1} & z_1 \\
    0 & 0 & 1
  \end{array}\right)
\end{equation} 
In other words: $T_{(z_1,t_1)}\psi(z_2,t_2,0)=\psi(z_1+z_2,t_1+t_2+2{\rm Im} \, (z_1 \cdot \overline{z_2}),0)$.

In the above terminology, the unipotent isometry (given by the projective action of) $T_{(z_1,t_1)}$ is 2-step (or vertical) if $z_1=0$ and 3-step (horizontal) otherwise. 

\smallskip

{\bf The hybridization construction:}

\smallskip

We will first embed the pair of Fuchsian groups into  $\SU(2,1)$ in the ball model of $H_\C^2$; there, two preferred orthogonal complex lines $L_1$ and $L_2$ are given by (the coordinate axes in the standard affine chart) $L_1=\pi({\rm Span}(e_1,e_3))$ and $L_2=\pi({\rm Span}(e_2,e_3))$, where $(e_1,e_2,e_3)$ denotes the canonical basis of $\C^3$ and $\pi:\C^3 - \{ 0\}\longrightarrow \C P^2$ the projectivization map. These intersect at the origin $O=\pi(e_3)$.
 
We will embed $\SU(1,1)$ in the stabilizer of each of these complex lines in the obvious block matrix form, namely via the injective homomorphisms:

\begin{equation}
  \begin{array}{llll}
    \iota_1: & \SU(1,1) & \longrightarrow &\SU(2,1) \\
    & \left(\begin{array}{cc}
      a & b \\
      c & d 
    \end{array}\right) &
    \longmapsto &
    \left(\begin{array}{ccc}
      a & 0 & b \\
      0 & 1 & 0 \\
      c & 0 & d
    \end{array}\right) 
  \end{array}
\end{equation}

\begin{equation}\label{rho2}
  \begin{array}{llll}
    \iota_2: & \SU(1,1) & \longrightarrow &\SU(2,1) \\
    & \left(\begin{array}{cc}
      a & b \\
      c & d 
    \end{array}\right) &
    \longmapsto &
    \left(\begin{array}{ccc}
      1 & 0 & 0 \\
      0 & a & b \\
      0 & c & d
    \end{array}\right) 
  \end{array}
\end{equation}

In the notation from the introduction, given two lattices $\Gamma_1,\Gamma_2$ in $\SU(1,1)$, we consider the hybrid $H(\Gamma_1,\Gamma_2)$ to be the projective image in $\PU(2,1)$ of $\langle \iota_1(\Gamma_1), \iota_2(\Gamma_2) \rangle < \SU(2,1)$.

\section{A hybrid subgroup of the Eisenstein-Picard modular group $\rm{PU}(2,1,\mathcal{O}_3)$}	

\begin{figure}[!hb]
    \centering
    \begin{tikzpicture}[x=1in,y=1in]
        \draw[thick] (0,0) circle (1);
        \coordinate (w) at ({1/2},{sqrt(3)/2});
        \coordinate (w-bar) at ({1/2},{-sqrt(3)/2});
        \draw[very thick] (w) arc (150:270:{1/sqrt(3)});
        \draw[very thick] (1,0) arc (90:210:{1/sqrt(3)});
        \draw[very thick] (0,0) -- (1,0);
        \draw[very thick] (0,0) -- (w);
        \draw[very thick] (0,0) -- (w-bar);
        \coordinate (R) at (0:0.25);
        \coordinate (U) at (0.5879,0.1823);
        \draw (R) node[anchor=north] {$R$};
        \draw (U) node[anchor=south west] {$U$};
        \draw[-latex] (R) arc (0:60:0.25);
        \draw[-latex] (U) arc (133.1736:226.8264:0.25);
    \end{tikzpicture}
    \caption{Fundamental domain for $\PU(1,1;\mathcal{O}_3)$}
\label{fig:d=3}
\end{figure}
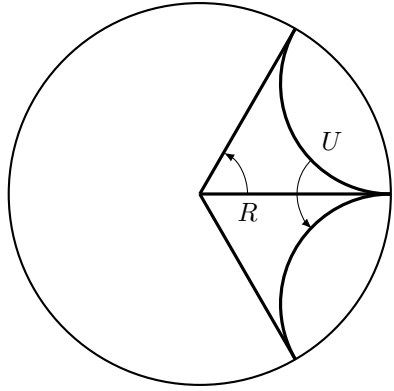

Using the method for finding generators from \cite{MP} (using 4 as \emph{covering depth}), we obtain the following generators for $\PU(1,1;\mathcal{O}_3)$ in the disk model of ${\rm H}_{\C}^1$:
\begin{align*}
R = \begin{pmatrix} -\omega^2 & 0 \\ 0 & 1 \end{pmatrix}, \qquad U = \begin{pmatrix} 1 + i\sqrt{3} & -i\sqrt{3} \\ i\sqrt{3} & 1 - i\sqrt{3} \end{pmatrix}, 
\end{align*}
denoting $\omega=\frac{-1+i\sqrt{3}}{2}$. Thus $\operatorname{U}(1,1;\mathcal{O}_3) = \langle -{\rm Id}, R, U \rangle$. The index-2 subgroup $\langle R^2, U \rangle < \PU(1,1;\mathcal{O}_3)$ lifts to $\SU(1,1;\mathcal{O}_3) = \langle -{\rm Id},E, U \rangle$ where
\begin{align*}
E=\omega^2 R^{-2} = \begin{pmatrix} \omega & 0 \\ 0 & \omega^2 \end{pmatrix}.
\end{align*}
Note that ${\rm PSU}(1,1;\mathcal{O}_3)=\langle E, U \rangle$ is a $(3,\infty,\infty)$ triangle group while $\PU(1,1;\mathcal{O}_3)$ is a $(2,6,\infty)$ triangle group. 

We consider the hybrid group $H(3)=H\left(\SU(1,1;\mathcal{O}_3),\SU(1,1;\mathcal{O}_3)\right)$, which by definition is generated by $\iota_j(E),\iota_j(U)$ and $\iota_j(-\rm{Id})$ for $j=1,2$. By the following observation it will suffice for our purposes to study the subgroup $\tilde{H}(3)=\langle \iota_1(E), \iota_2(E), \iota_1(U), \iota_2(U) \rangle$, which is a hybrid of two copies of ${\rm PSU}(1,1;\mathcal{O}_3)$ (so, of two $(3,\infty,\infty)$ triangle groups, ignoring the central extension by $-{\rm Id}$). 

\begin{lem}
  \label{HtildeinH} $\tilde{H}(3)$ is normal in $H(3)$, with index dividing 4.
\end{lem}

\Pf Since the $\iota_j(E), \iota_j(-\rm{Id})$ are diagonal they all commute; moreover for each $j=1,2$ $\iota_j(-\rm{Id})$ commutes with $\iota_j(E), \iota_j(U)$. It suffices therefore to show that $\iota_j(-Id)$ conjugates $\iota_{3-j}(U)$ into  $\tilde{H}(3)$ for $j=1,2$; a straightforward computation gives: $\iota_j(-{\rm Id}) \iota_{3-j}(U) \iota_j(-{\rm Id}) =\iota_{3-j}((EUE)^{-1})$. Therefore $\tilde{H}(3)$ is normal in $H(3)$; the quotient group is a quotient of $\Z/2\Z \times \Z/2\Z$ as it is generated by the images of $\iota_1(-\rm{Id})$ and $\iota_2(-\rm{Id})$ which have order 2 and commute. \EPf

It will be more convenient for us to work in the Siegel model, in other words to conjugate by the Cayley transform $J$ given in (\ref{cayley}); we will abuse notation slightly by denoting again $H(3), \tilde{H}(3)$ the conjugates $J^{-1}H(3)J, J^{-1}\tilde{H}(3)J$. Concretely, we consider $\tilde{H}(3)= \langle E_1,U_1,E_2,U_2 \rangle$, with:
\begin{align*}
  E_1 = J^{-1} \iota_1(E) J
  &= 
  \begin{pmatrix}
    \omega^2 & \omega^2-1 & \omega+2\\
    i \sqrt{3} & 1+ i \sqrt{3} & \omega^2-1 \\
    i \sqrt{3} & i \sqrt{3} & \omega^2
  \end{pmatrix},
  &
  U_1 = J^{-1} \iota_1(U) J
  &= 
  \begin{pmatrix}
    1 & 0 & i \sqrt{3} \\
    0 & 1 & 0 \\
    0 & 0 & 1 
  \end{pmatrix},
  \\
  E_2 = J^{-1} \iota_2(E) J
  &=
  \begin{pmatrix}
    \omega^2 & -i \sqrt{3} & i \sqrt{3} \\
   \omega+2& 1+ i \sqrt{3} & -i \sqrt{3} \\
   \omega+2 &\omega+2 & \omega^2
  \end{pmatrix},
  &
  U_2 = J^{-1} \iota_2(U) J
  &=
  \begin{pmatrix}
    1 & 0 & 0 \\
    0 & 1 & 0 \\
    i \sqrt{3} & 0 & 1
  \end{pmatrix}.
\end{align*}
Likewise we will denote $I_1 = J^{-1} \iota_{1}(-{\rm Id}) J$ and $I_2 = J^{-1} \iota_2(-{\rm Id}) J$.

{\bf Remark:} Since $E_1, E_2$ are both regular elliptic of order 3 with the same eigenspaces, they are either equal or inverse of each other. It turns out that $E_2=E_1^{-1}$ in $\PU(2,1)$ (the matrices satisfy $E_2=\omega E_1^{-1}$). We will therefore omit the generator $E_2$ from now on.

In \cite{FP} the authors determine that the Eisenstein-Picard modular group $\PU(2,1;\mathcal{O}_3)$ has presentation:

\begin{align*}
  \PU(2,1;\mathcal{O}_3)=\left\langle P,Q,R \;|\;R^2,\, (QP^{-1})^{6},\,PQ^{-1}RQP^{-1}R,\,P^{3}Q^{-2},\,(RP)^{3}\right\rangle , \ \ {\rm where}
\end{align*}

\begin{align*}
  P=\begin{pmatrix}
  1 & 1 &\omega\\
  0 & \omega &-\omega \\
  0 & 0 & 1
  \end{pmatrix},
  \qquad
  Q=\begin{pmatrix}
  1 & 1 & \omega \\
  0 & -1 & 1 \\
  0 & 0 & 1
  \end{pmatrix},
  \qquad
  R=\begin{pmatrix}
  0 & 0 & 1 \\
  0 & -1 & 0 \\
  1 & 0 & 0
  \end{pmatrix}.
\end{align*}

The following result can be checked by straightforward matrix computation:

\begin{lem}\label{gens3}
  The generators for the hybrid $\tilde{H}(3)$ can be expressed in terms of the Falbel-Parker generators for $\PU(2,1;\mathcal{O}_3)$ as follows:
  \begin{align*}
    U_1 &= Q^2, \\
    U_2 &= R Q^2 R, \\
    E_1 &= P^2 (R Q^2)^2 P^{-2}.
  \end{align*}    
\end{lem}

\begin{lem}\label{normal3}
  The hybrid $\tilde{H}(3)$ is a normal subgroup of $\PU(2,1;\mathcal{O}_3)$.
\end{lem}

\Pf It suffices to check that generators of $\PU(2,1;\mathcal{O}_3)$ conjugate generators of $\tilde{H}(3)$ into $\tilde{H}(3)$. Straightforward computations give the following:

$$\begin{array}{lcr}
  \begin{array}{rcl}
    P^{-1}U_1 P &= &U_1 \\
    Q^{-1}U_1 Q &= & U_1 \\
    R^{-1}U_1 R &= &U_2 
  \end{array}
  & 
  \begin{array}{rcl}
    P^{-1}U_2 P &= &U_1^{-1}E_1 \\
    Q^{-1}U_2 Q &= &U_1^{-1}E_1\\
    R^{-1}U_2 R &= &U_1 
  \end{array}
  &
  \begin{array}{rcl}
    P^{-1}E_1 P &= &U_2^{-1} E_1^{-1} U_1 \\
    Q^{-1}E_1 Q &= &U_2 U_1 \\
    R^{-1}E_1 R &= &E_1^{-1}
  \end{array} 
\end{array}$$
\EPf

We can then form the quotient group $G=\PU(2,1;\mathcal{O}_3)/ \tilde{H}(3)$, which by Lemma~\ref{gens3} has presentation:

\begin{align*}
  G=\PU(2,1;\mathcal{O}_3)/ \tilde{H}(3)  
  = \left\langle P, Q, R \;\bigg|\; 
  \begin{aligned}
    &R^2, (QP^{-1})^6, PQ^{-1}RQP^{-1}R, \\
    &P^3Q^{-2}, (RP)^3, Q^2
  \end{aligned}
  \right\rangle.
\end{align*}
(Note that the relation $Q^{2}$ makes the other three relators corresponding to the generators of $H(3)$ superfluous). The Tietze transformation $a=PQ^{-1}$, $b=Q$, $c=R$, yields the following presentation for $G$:

\begin{align*}
  G = \left\langle a,b,c \;|\; c^2,\, a^6,\, [a,c],\, (ab)^3,\,(cab)^3,\, b^2\right\rangle.
\end{align*}

\begin{thm}\label{d=3main}
  The hybrid $H(3)$ has infinite index in $\PU(2,1,\mathcal{O}_3)$.
\end{thm}

Note that $G / \langle\!\langle c \rangle\!\rangle$ is the $(2,3,6)$ triangle group, hence infinite. Therefore $G$ is also infinite, in other words $\tilde{H}(3)$ has infinite index in $\PU(2,1,\mathcal{O}_3)$. By Lemma~\ref{HtildeinH}, $H(3)$ also has infinite index in $\PU(2,1,\mathcal{O}_3)$. \EPf

\begin{cor}\label{d=3Zdense}
  The hybrid $H(3)$ is a thin subgroup of $\PU(2,1,\mathcal{O}_3)$.
\end{cor}

\Pf The only additional statement is that $H(3)$ is Zariski-dense in $\PU(2,1)$, which is simple to see in rank 1, as it reduces essentially to irreducibility. Indeed, by \cite{CG} if a discrete subgroup $\Gamma$ is not Zariski-dense then it preserves a strict subspace of ${\rm H}^2_\C$ or it fixes a point on $\partial_\infty {\rm H}^2_\C$. This is easily seen not to be the case, as $E_1$ does not preserve the unique complex line preserved by both $U_1$ and $U_2$. (This will also follow from the fact that $H(3)$ has full limit set, see Section 6).  \EPf

We conclude this section with a few remarks about the algebraic structure of the hybrid $H(3)$. We do not know a complete presentation for $H(3)$, in fact it may be non-finitely presented as far as we know (see \cite{K2} and Proposition~4.2 of \cite{FW}). The following observations are obtained by direct computation using the generators in matrix form.

\begin{lem}\label{rels3}
  The following relations hold between the generators $E_1,U_1,U_2$ for $H(3)$: 
  $$E_1^3=(U_1U_2)^3=(E_1U_1^{-1}U_2)^3=(E_1U_2U_1^{-1})^3=(E_1^{-1}U_1U_2^{-1})^3=(E_1^{-1}U_2^{-1}U_1)^3=1.$$
\end{lem}

\begin{cor}\label{notamalgam}
  The hybrid $H(3)$ has finite abelianization; in particular it is not isomorphic to the amalgamated product of $\iota_1(\SU(1,1,\mathcal{O}_3))$ and $\iota_2(\SU(1,1,\mathcal{O}_3))$ over their (projective) intersection.
\end{cor}

\Pf Recall that $H(3)$ is generated by $E_1, U_1, U_2, I_1, I_2$ with $I_1, I_2$ both of order $2$. Observe that by Lemma~\ref{rels3}, the following relations hold in the abelianization $H(3)^{ab}$ (we slightly abuse notation by using the same symbol for elements of $H(3)$ and their image in $H(3)^{ab}$): $ E_1^3=U_1^6=1, \ U_1^3=U_2^3$. Therefore $H(3)^{ab}$ is a quotient of $(\Z/2\Z)^2 \times \Z/3\Z \times (\Z/6\Z)^2$. 

For the second statement, recall from the beginning of Section 3 that $\SU(1,1,\mathcal{O}_3)$ is a central $\Z/2\Z$-extension of ${\rm PSU}(1,1,\mathcal{O}_3)$, a $(3,\infty,\infty)$ triangle group. Therefore $\SU(1,1,\mathcal{O}_3)$ has abelianization $\Z \times \Z/6\Z$. Now the projective images of $\iota_1(\SU(1,1,\mathcal{O}_3))$ and $\iota_2(\SU(1,1,\mathcal{O}_3))$ intersect along the cyclic group of order 3 generated by $E_1$ (indeed, this intersection is contained in the diagonal subgroup of ${\rm PU}(2,1,\mathcal{O}_3)$, and $E_1,E_1^2$ are the only non-trivial diagonal elements that can be scaled to $\{ 1,u,u^{-1}\}$ for some $u \in {\rm U}(1)$). Therefore the amalgamated product has abelianization $\Z^2 \times \Z/3\Z \times \Z/6\Z$. \EPf

It is interesting to note that this also tells us the behavior of a related hybrid group, namely the hybrid of two 
$(2,6,\infty)$ triangle groups, rather than $(3,\infty,\infty)$ (which each $(2,6,\infty)$ group contains with index 2). A simple way to view this new hybrid $H'(3)$ as a subgroup of $\Gamma(3)=\PU(2,1,\mathcal{O}_3)$ containing the previous hybrid $\tilde{H(3)}$ is to take the obvious square root of the previous generator $E_1$ in terms of the Falbel-Parker generators, in other words to take $H'(3)$ to be generated by $E_1'=P^2(RQ^2)P^{-2}$, and $U_1=Q^2,U_2=RQ^2R$ unchanged.

\begin{lem}\label{comm3}
  The hybrid $H'(3)$ is contained in $[\Gamma(3),\Gamma(3)]$.
\end{lem}

\Pf From the Falbel-Parker presentation for $\Gamma_3$ we get (abusing notation slightly again by using the same symbol for elements of $\Gamma(3)$ and their image in $\Gamma(3)^{ab})$:
$$\Gamma(3)^{ab}=\Gamma(3)/[\Gamma(3),\Gamma(3)] = \langle P,Q,R \ | \ R=P^3=Q^2=[P,Q]=1  \rangle.$$
The result then follows by noting that the generators listed above for $H'(3)$ all become trivial in the abelianization. \EPf

The following is Lemma~6 of \cite{FW}.
\begin{lem}\label{comm3ab}
  The commutator subgroup $[\Gamma(3),\Gamma(3)]$ has abelianization $\Z \oplus \Z$.
\end{lem}

\begin{lem}\label{H'fab}
  The hybrid $H'(3)$ has finite abelianization.
\end{lem}

\Pf This follows from the relations given in Lemma~\ref{rels3} by noting that $H'(3)$ is generated by $E'_1,U_1,U_2$ with $(E'_1)^2=E_1$. \EPf

The following is well known but we include it for completeness:

\begin{lem}\label{fiabs}
  If $K_1 < K_2$ are two groups with $[K_2:K_1]$ and $K_1^{ab}$ finite, then $K_2^{ab}$ is finite. 
\end{lem}

\Pf Denote $i$ the inclusion map from $K_1$ into $K_2$, and $\pi_i:K_i\longrightarrow K_i^{ab}$ the quotient maps for $i=1,2$. Then $\pi_2 \circ i$ is a homomorphism from $K_1$ to an abelian group, so by the universal property of abelianizations $\pi_2 \circ i$ factors through $K_1^{ab}$, i.e. there is a homomorphism $i_*:K_1^{ab} \longrightarrow K_2^{ab}$ such that $i_* \circ \pi_1=\pi_2 \circ i$. Since $K_1=i(K_1)$ has finite index in $K_2$ by assumption and $\pi_2$ is surjective, $\pi_2(K_1)=i_*(\pi_1(K_1))=i_*(K_1^{ab})$ has finite index in $K_2^{ab}$. The result follows since $K_1^{ab}$ is finite. \EPf

Combining Lemmas~\ref{comm3}, \ref{comm3ab}, \ref{H'fab} and \ref{fiabs} gives the following:

\begin{cor}\label{d=3other}
  The hybrid $H'(3)$ has infinite index in $[\Gamma(3),\Gamma(3)]$, hence also in $\Gamma(3)$.
\end{cor}

It is interesting to note that, in contrast with the previous hybrid $H(3)$ which was normal in $\Gamma(3)$, $H'(3)$ now has infinite index in its normal closure $\langle \langle H'(3) \rangle \rangle =\Gamma(3)$ in $\Gamma(3)$ (the presentation of $\Gamma(3)/\langle \langle H'(3) \rangle \rangle$ obtained by adding the generators of $H'(3)$ to the presentation for $\Gamma(3)$ now gives the trivial group).

\section{A hybrid subgroup of the Gauss-Picard modular group $\rm{PU}(2,1,\mathcal{O}_1)$}

\begin{figure}[!hb]
    \centering
    \begin{tikzpicture}[x=1in,y=1in]
        \draw[thick] (0,0) circle (1);
        \coordinate (i) at (0,1);
        \coordinate (i-bar) at (0,-1);
        \draw[very thick] (i) arc (180:270:1);
        \draw[very thick] (1,0) arc (90:180:1);
        \draw[very thick] (0,0) -- (1,0);
        \draw[very thick] (0,0) -- (i);
        \draw[very thick] (0,0) -- (i-bar);
        \coordinate (R) at (0:0.20);
        \coordinate (U) at (0.5294,0.1176);
        \draw (R) node[anchor=north] {$R$};
        \draw (U) node[anchor=south west] {$U$};
        \draw[-latex] (R) arc (0:90:0.20);
        \draw[-latex] (U) arc (151.9275:208.0725:0.25);
    \end{tikzpicture}
\caption{Fundamental domain for $\PU(1,1;\mathcal{O}_1)$}
\label{fig:d=1}
\end{figure}
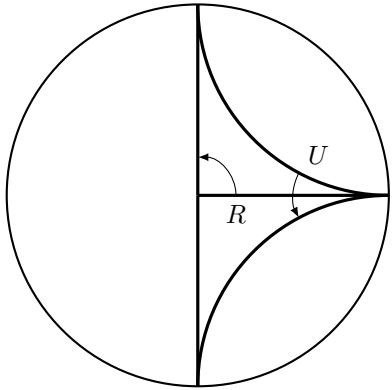

Using the method for finding generators from \cite{MP} (using again 4 as \emph{covering depth}), we obtain the following generators for $\PU(1,1;\mathcal{O}_1)$ in the disk model of $H_C^1$:
\begin{align*}
R = \begin{pmatrix} i & 0 \\ 0 & 1 \end{pmatrix}, \qquad U = \begin{pmatrix} 1 + i & -i \\ i & 1 - i \end{pmatrix}.
\end{align*}
Thus $\operatorname{U}(1,1;\mathcal{O}_1) = \langle -{\rm Id}, R, U \rangle$. The index-2 subgroup $\langle R^2, U \rangle < \PU(1,1;\mathcal{O}_1)$ lifts to $\SU(1,1;\mathcal{O}_1)$, which is therefore generated by $-{\rm Id},E,$ and $T$, with:
\begin{align*}
E=i R^{-2} = \begin{pmatrix} -i & 0 \\ 0 & i \end{pmatrix}.
\end{align*}
Since $E^{2} = -{\rm Id}$, the group $\SU(1,1;\mathcal{O}_1)$ is generated by $E$ and $U$.

Note that $\SU(1,1;\mathcal{O}_1)$ is a central $\Z/2\Z$-extension of a $(2,\infty,\infty)$ triangle group while $\PU(1,1;\mathcal{O}_1)$ is a $(2,4,\infty)$ triangle group.

We now consider the hybrid group $H\left(\SU(1,1;\mathcal{O}_1),\SU(1,1;\mathcal{O}_1)\right)$, which by definition is generated by  $\iota_1(E)$, $\iota_1(U)$, $\iota_2(E)$ and $\iota_2(U)$. It will be again more convenient for us to work in the Siegel model, in other words to conjugate by the Cayley transform $J$ given in (\ref{cayley}). We thus consider the group $H(1)= \langle E_1,U_1,E_2,U_2 \rangle$, where:

\begin{align*}
  E_1 = J^{-1}\iota_1(E)J
  &=
  \begin{pmatrix}
    i & -1+i & 1-i \\
    -2i & 1-2i & -1+i \\
    -2i & -2i & i
  \end{pmatrix},
  &
  U_1 =J^{-1}\iota_1(U)J
  &=
  \begin{pmatrix}
    1 & 0 & i \\
    0 & 1 & 0 \\
    0 & 0 & 1
  \end{pmatrix}, \\
  E_2 = J^{-1}\iota_2(E)J
  &=
  \begin{pmatrix}
    i & 2i & -2i \\
    1-i & 1-2i & 2i \\
    1-i & 1-i & i
  \end{pmatrix},
  &
  U_2 = J^{-1}\iota_2(U)J
  &=
  \begin{pmatrix}
    1 & 0 & 0 \\
    0 & 1 & 0 \\
    i & 0 & 1
  \end{pmatrix}.
\end{align*}

The following presentation for the Gauss-Picard lattice $\PU(2,1;\mathcal{O}_1)$ was first given in \cite{FFP} (Theorem 4):

\begin{align*}
  \PU(2,1;\mathcal{O}_1) = \left\langle I_0, Q, T \; \bigg| \;
  \begin{aligned}
    I_0^2,\, Q^2,\, (I_0 Q)^3,\, (I_0 T)^{12},\, (I_0 QT)^{8},\, [(I_0 T)^3,T],\, [Q,T]
  \end{aligned}
  \right\rangle
\end{align*}

where:

\begin{align*}
  I_0 = \begin{pmatrix} 0 & 0 & 1 \\ 0 & -1 & 0 \\ 1 & 0 & 0 \end{pmatrix},
  \qquad
  Q = \begin{pmatrix} 1 & 1-i & -1 \\ 0 & -1 & 1 + i \\ 0 & 0 & 1 \end{pmatrix},
  \qquad
  T = \begin{pmatrix} 1 & 0 & i \\ 0 & 1 & 0 \\ 0 & 0 & 1 \end{pmatrix}.
\end{align*}



The following result can be checked by straightforward matrix computation; the expression for $E_1$ was found by a computer search, using a different generating set.

\begin{lem}\label{d-1gens}
  The generators for the hybrid $H(1)$ can be expressed in terms of the above generators for $\PU(2,1;\mathcal{O}_1)$ as follows:
  \begin{align*}
    U_1 &= T, \\
    U_2 &= I_0 U_2 I_0, \\
    E_1 &= T^{-1} Q (I_0 T)^3 I_0 (T (I_0 T)^{-3} Q)^2 I_0, \\
    E_2 &= I_0 E_1 I_0.
  \end{align*}
\end{lem}

\begin{lem}\label{normal1}
  The hybrid $H(1)$ is a normal subgroup of $\PU(2,1;\mathcal{O}_1)$.
\end{lem}

\Pf Since $U_1 = T$ and the generators come in pairs conjugate by $I_0$, it suffices to check conjugation by $Q$:
\begin{align*}
  Q^{-1} U_1 Q &= U_1 \\
  Q^{-1} U_2 Q &= (U_1 E_1) U_2 (U_1 E_1)^{-1} \\
  Q^{-1} E_1 Q &= (U_2 U_1) E_2 (U_2 U_1)^{-1} \\
  Q^{-1} E_2 Q &= (U_2 U_1) E_1 (U_2 U_1)^{-1}
\end{align*}
\EPf

\begin{thm}\label{d=1main}
  The hybrid $H(1)$ has index 2 in the full Gauss-Picard lattice $\PU(2,1;\mathcal{O}_1)$.
\end{thm}

\Pf Since $H(1)$ is normal in $\PU(2,1;\mathcal{O}_1)$, a presentation for the quotient $\PU(2,1;\mathcal{O}_1)/H(1)$ is obtained from the presentation for $\PU(2,1;\mathcal{O}_1)$, to which we add as relations the generators of the subgroup $H(1)$ written as the words given in Lemma~\ref{d-1gens}.
\begin{align*}
  \PU(2,1;\mathcal{O}_1)/H(1)
  &= \left\langle I_0, Q, T \; \bigg| \;
  \begin{aligned}
    I_0^2,\, Q^2,\, (I_0 Q)^3,\, (I_0 T)^{12},\, (I_0 QT)^{8},\, [(I_0 T)^3,T],\, [Q,T],\, T,\, I_0 T I_0, \\
    T^{-1} Q (I_0 T)^3 I_0 (T (I_0 T)^{-3} Q)^2 I_0,\, I_0 T^{-1} Q (I_0 T)^3 I_0 (T (I_0 T)^{-3} Q)^2
  \end{aligned}
  \right\rangle
\end{align*}
Since $T = 1$ in the quotient, the relations $(I_0 Q)^3 = 1$ and $(I_0 Q T)^8 = 1$ imply that $I_0 = Q$. The remaining relations are then satisfied, thus the presentation for the quotient simplifies to
\begin{align*}
  \PU(2,1;\mathcal{O}_1)/H(1) = \left\langle I_0, Q, T \;|\; I_0 = Q, I_{0}^2 = T = 1 \right\rangle = \Z/2\Z .
\end{align*}
\EPf

We now consider the related hybrid $H'(1)$ as in the case of $d=3$, namely taking $H'(1)$ to be the hybrid of two copies of ($\Z/2\Z$-extensions of) the Fuchsian triangle group $(2,4,\infty)$, rather than $(2,\infty,\infty)$. We immediately get the following result by noting that $H'(1)$ contains $H(1)$, which has index 2 in the full lattice $\Gamma(1)$, as well as a new element of order 4 not belonging to $H(1)$.

\begin{cor}\label{d=1other}
  The hybrid $H'(1)$ is equal to the full lattice $\Gamma(1)=\PU(2,1;\mathcal{O}_1)$.
\end{cor}

\section{A hybrid subgroup of the Picard modular group $\rm{PU}(2,1,\mathcal{O}_7)$}

\begin{figure}[!hb]
    \centering
    \begin{tikzpicture}[x=1in,y=1in]
        \draw[thick] (0,0) circle (1);
        \begin{scope}
        \coordinate (P1) at (0,0.4569);
        \coordinate (P2) at (0.625,0.3307);
        \coordinate (P3) at (1,0);
        \coordinate (P4) at (0.625,-0.3307);
        \coordinate (P5) at (0,-0.4569);
        \draw[very thick] (P1) arc (240 : 277.1808 : 1);
        \draw[very thick] (P2) arc (187.1807 : 270 : 0.37796);
        \draw[very thick] (0,0) -- (P3);
        \draw[very thick] (P3) arc (90: 172.8193 : 0.37796);
        \draw[very thick] (P4) arc (82.8192 : 120 : 1); 
        \fill (0,0) circle (2pt);
        \coordinate (B) at (0:0.15);
        \coordinate (U) at (0.75,0.09448);
        \coordinate (A) at (0.4251,0.3257);
        \draw (B) node[anchor=north] {$B$};
        \draw (U) node[anchor=south west] {$U$};
        \draw (A) node[anchor=south] {$A$};
        \draw[-latex] (B) arc (0:180:0.15);
        \draw[-latex] (U) arc (138.592 : 221.408 : 0.14285);
        \draw[-latex] (A) arc (175.697 : 184.303 : 4.3411);
        \end{scope}
        \begin{scope}[xscale=-1,yscale=1]
        \coordinate (P1) at (0,0.4569);
        \coordinate (P2) at (0.625,0.3307);
        \coordinate (P3) at (1,0);
        \coordinate (P4) at (0.625,-0.3307);
        \coordinate (P5) at (0,-0.4569);
        \draw[very thick] (P1) arc (240 : 277.1808 : 1);
        \draw[very thick] (P2) arc (187.1807 : 270 : 0.37796);
        \draw[very thick] (0,0) -- (P3);
        \draw[very thick] (P3) arc (90: 172.8193 : 0.37796);
        \draw[very thick] (P4) arc (82.8192 : 120 : 1); 
        \end{scope}
    \end{tikzpicture}
\caption{Fundamental domain for $\PU(1,1;\mathcal{O}_7)$}
\label{fig:d=7}
\end{figure}
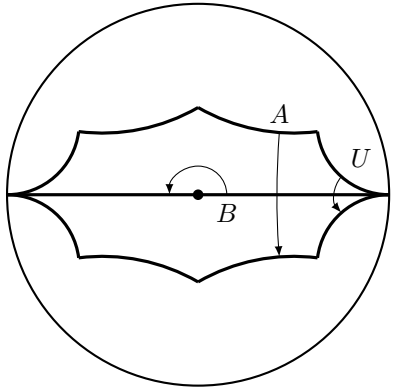

Using the method for finding generators from \cite{MP} (using 9 as \emph{covering depth}), we obtain the following generators for $\PU(1,1;\mathcal{O}_7)$ in the disk model of ${\rm H}_{\C}^1$:
\begin{align*}
  U = \begin{pmatrix} 1 + i \sqrt{7}& -i\sqrt{7} \\ i\sqrt{7} & 1-i\sqrt{7} \end{pmatrix}, \qquad 
  A = \begin{pmatrix} -\frac{1}{2} + i\frac{\sqrt{7}}{2} & 1 \\ -1 & \frac{1}{2} + i\frac{\sqrt{7}}{2} \end{pmatrix},  \qquad
  B = \begin{pmatrix} -1 & 0 \\ 0 & 1 \end{pmatrix}.
\end{align*}
Thus $\operatorname{U}(1,1;\mathcal{O}_7) = \langle -\operatorname{Id},U,A,B \rangle$. As previously, the index-2 subgroup $\langle U, BA, AB, BUB \rangle < \PU(1,1;\mathcal{O}_7)$ lifts to $\SU(1,1;\mathcal{O}_7)$, which is therefore generated by $ -\operatorname{Id}, U, BA, AB,$ and $BUB$. 

When embedding $\SU(1,1;\mathcal{O}_7)$ into $\PU(2,1)$ via $\iota_1$ and $\iota_2$, it turns out that $\iota_1(-{\rm Id}) = \iota_2(B)$ and $\iota_2(-\rm{Id}) = \iota_1(B)$, and so the hybrid group $H(7) = H(\SU(1,1;\mathcal{O}_7),\SU(1,1;\mathcal{O}_7))$ admits the following simplified set of generators:

\begin{align*}
  U_1 = J^{-1}\iota_1(U)J
   &=
  \begin{pmatrix}
    1 & 0 & i \sqrt{7} \\
    0 & 1 & 0 \\
    0 & 0 & 1
  \end{pmatrix},
  &
  U_2 = J^{-1}\iota_2(U)J
  &=
  \begin{pmatrix}
    1 & 0 & 0 \\
    0 & 1 & 0 \\
    i \sqrt{7} & 0 & 1
  \end{pmatrix},
  \\
  A_1 = J^{-1}\iota_1(A)J
  &=
  \begin{pmatrix}
    -\frac{1}{2}+i\frac{\sqrt{7}}{2} & -\frac{3}{2}+i\frac{\sqrt{7}}{2} & \frac{1}{2}-i\frac{\sqrt{7}}{2} \\
    1 & 2 & -\frac{3}{2}+i\frac{\sqrt{7}}{2} \\
    1 & 1 & -\frac{1}{2}+i\frac{\sqrt{7}}{2} \\
  \end{pmatrix},
  &
  B_1 = J^{-1}\iota_1(B)J
  &=
  \begin{pmatrix}
    1 & 0 & 0 \\
    -2 & -1 & 0 \\
    -2 & -2 & 1
  \end{pmatrix},
  \\
  A_2 = J^{-1}\iota_2(A)J
  &=
  \begin{pmatrix}
    -\frac{1}{2} + i\frac{\sqrt{7}}{2} & -1 & 1 \\
    \frac{3}{2} - i\frac{\sqrt{7}}{2} & 2 & -1 \\
    \frac{1}{2} - i\frac{\sqrt{7}}{2} & \frac{3}{2}-\frac{i \sqrt{7}}{2} & -\frac{1}{2}+\frac{i \sqrt{7}}{2} \\
  \end{pmatrix},
  &
  B_2 = J^{-1}\iota_2(U)J
  &=
  \begin{pmatrix}
    1 & 2 & -2 \\
    0 & -1 & 2 \\
    0 & 0 & 1
  \end{pmatrix}.
\end{align*}

In ~\cite{MP} the authors determine that $\PU(2,1;\mathcal{O}_7)$ has presentation
\begin{align}\label{MPpres}
  \PU(2,1;\mathcal{O}_7)= \left\langle T_1, R, I\; \bigg| \;
  \begin{aligned}
    R^2 =  I^2 = (RI)^2 = R  T_1^{-1}  R  T_1 RT_1RT_1^{-1} = 
    (T_1  I  T_1^{-1}  R)^4 = (T_1^{-1}  I  T_1  R)^4 =   \\
    T_1^{-1}  I  T_1^{-1}  I  T_1  I  T_1  I  T_1^{-3}  I  T_1 
    I  T_1  I  T_1^{-1}  I  T_1^{-1} = 
    (T_1^{-1}  I  T_1  I  T_1  I  T_1^{-1}  I  T_1^{-1}  I)^2 = \\
     (I  T_1^{-1}  R)^7 =
    T_1^{-1}  I  T_1  I  T_1  I  T_1^{-2}  I  T_1^{-1}  I  T_1 
    I  T_1^2  I  T_1^{-1}  I  T_1^{-1}  I  T_1  I =  \\
    T_1^{-1}  I  T_1  I  T_1  I  R  T_1  I  R  T_1  I 
    T_1  I  T_1^{-1}  I  T_1^{-1}  I  T_1  R  T_1^{-1}  I  R 
    T_1^{-1}  I = \\
    R  T_1  I  R  T_1  I  T_1  I  T_1^{-1}  I  T_1^{-1}  I
     R  T_1^{-1}  I  R  T_1^{-1}  I  T_1^{-1}  I  T_1  I  T_1 
    I  T_1^{-1} = \\
    R  T_1  I  R  T_1  R  T_1^{-1}  I  T_1  I  T_1  I 
    R  T_1  I  T_1  I  T_1^{-1}  R  T_1  R  I  T_1  R 
    T_1^{-1}  I  T_1  I  T_1  I  T_1  I  T_1^{-1} = 1
 \end{aligned}
  \right\rangle
\end{align}
where

\begin{align*}
  T_1 &= \begin{pmatrix} 
    1 & -1 & -\frac{1}{2} + i\frac{\sqrt{7}}{2} \\ 
    0 & 1 & 1 \\ 
    0 & 0 & 1 \end{pmatrix}, &
  R &= \begin{pmatrix}
    1 & 0 & 0 \\
    0 & -1 & 0 \\
    0 & 0 & 1 \end{pmatrix}, &
  I &= \begin{pmatrix}
    0 & 0 & 1 \\
    0 & -1 & 0 \\
    1 & 0 & 0 \end{pmatrix}. 
\end{align*}

%
%
The following result can be checked by straightforward matrix computation: 
\begin{lem}\label{d=7gens}
In terms of the above generators, the generators for $H(7)$ can be written as follows:
\begin{align*}
  U_1 &= (R T_1)^2, \\
  U_2 &= I U_1 I, \\
  A_1 &= T_1 I T_1 R, \\
  A_2 &= I A_1 I,\\
  B_1 &= (IT_1) R (I T_1)^{-1},\\
  B_2 &= I B_1 I.
\end{align*}
\end{lem}

\begin{lem}\label{d=7normal}
  The hybrid $H(7)$ is a normal subgroup of $\PU(2,1;\mathcal{O}_7)$.
\end{lem}

\Pf
We have that $R = (A_1 A_2 B_1 A_1 B_2)^{-1} B_1 (A_1 A_2 B_1 A_1 B_2) \in H(7)$ and since the generators of $H(7)$ come in pairs conjugate by $I$, it suffices to check conjugation by $T_1$:
\begin{align*}
  T_1^{-1}A_1 T_1 &= (A_1 A_2^{-1} B_2 A_2^{-1} A_1)^{-1} A_2 (A_1 A_2^{-1} B_2 A_2^{-1} A_1) \\
  T_1^{-1}A_2 T_1 &= (B_2 A_2 A_1^{-1} A_2^{-1} B_1)^{-1} A_2 (B_2 A_2 A_1^{-1} A_2^{-1} B_1) \\
  T_1^{-1}B_1 T_1 &= (A_1^{-1} A_2^{-1} B_1)^{-1} B_1 (A_1^{-1} A_2^{-1} B_1) \\
  T_1^{-1}B_2 T_1 &= R \\
  T_1^{-1}U_1 T_1 &= U_1 \\
  T_1^{-1}U_2 T_1 &= (A_1^{2} A_2^{-1} B_2 A_2^{-1} A_1)^{-1} U_2 (A_1^{2} A_2^{-1} B_2 A_2^{-1} A_1) \\
\end{align*}
\EPf

\begin{thm}\label{d=7main}
  The hybrid ${\rm H}(7)$ is the full lattice $\PU(2,1;\mathcal{O}_7)$.
\end{thm}

\Pf As before, since $H(7)$ is normal in $\PU(2,1;\mathcal{O}_7)$, we obtain a presentation for the quotient $\PU(2,1;\mathcal{O}_7)/H(7)$ by adding the words from Lemma~\ref{d=7gens} as relations to the presentation (\ref{MPpres}). As noted in the proof of Lemma~\ref{d=7normal}, $R \in H(7)$ so $R=1$ in the quotient. As well, the relation corresponding to $A_1$ implies that $I = T_1^2$ in the quotient. The following relations in the presentation (\ref{MPpres})
\begin{align*}
I^2 = T_1^{-1} I T_1^{-1} I T_1 I T_1 I T_1^{-3} I T_1 I T_1 I T_1^{-1} I T_1^{-1} = 1
\end{align*}
imply that $T_1^4 = T_1^{13} = 1$ in the quotient, hence $T_1 = 1$. Thus
$\PU(2,1;\mathcal{O}_7)/H(7)$ is trivial. \EPf

\section{Limit sets and geometrical finiteness}

\subsection{Limit sets}

We first briefly recall the definition and two classical facts about limit sets of discrete groups of isometries, see e.g. \cite{B2} or \cite{K3}. The space we consider in this paper is the complex hyperbolic plane ${\rm H}^2_\C$, but these definitions and facts hold more generally in any negatively curved symmetric space (so, hyperbolic space of any dimension over the real or complex numbers or quaternions, or hyperbolic plane over the octonions).

{\bf Definiton:} Let $X$ be a negatively curved symmetric space, $\partial_\infty X$ its boundary at infinity, and $\Gamma$ a discrete subgroup of ${\rm Isom}(X)$. The \emph{limit set} $\Lambda(\Gamma)$ of $\Gamma$ is defined as the set of accumulation points in $\partial_\infty X$ of the orbit $\Gamma x_0$ for any choice of $x_0 \in X$; this does not depend on the choice of $x_0$.

A discrete subgroup $\Gamma$ of ${\rm Isom}(X)$ is called \emph{elementary} if $\Gamma$ has a global fixed point in $X \cup \partial_{\infty} X$ or preserves a geodesic. When $\Gamma$ is non-elementary, the limit set $\Lambda(\Gamma)$ is the minimal (nonempty) closed $\Gamma$-invariant subset of $\partial_\infty X$, in fact the orbit $\Gamma p_\infty$ is dense in $\Lambda(\Gamma)$ for any $p_\infty \in \Lambda(\Gamma)$. We will use the following two classical properties of limit sets:

\begin{prop} Let $X$ be a negatively curved symmetric space and $\Gamma$ a discrete subgroup of ${\rm Isom}(X)$. 
  \begin{itemize}
  \item[(a)] If $\Gamma$ is a lattice in ${\rm Isom}(X)$ then $\Lambda(\Gamma)=\partial_\infty X$.
  \item[(b)] If $\Gamma'$ is a nonelementary normal subgroup of $\Gamma$ then $\Lambda(\Gamma')=\Lambda(\Gamma)$.
  \end{itemize}
\end{prop}

The following result is an immediate consequence of this and Lemmas~\ref{normal3}, \ref{normal1} (or Theorem~\ref{d=1main}).

\begin{prop}\label{limitsets} For $d=1,3,7$ the hybrid $H(d)$ has full limit set: $\Lambda(H(d))=\partial_\infty {\rm H^2_\C} \simeq S^3$.
\end{prop}

\subsection{Geometrical finiteness}

The original notion of geometrical finiteness for a Kleinian group $\Gamma < {\rm Isom \, (H}^3_\R)$ was to admit a finite-sided polyhedral fundamental domain. This was later shown to admit several equivalent formulations, then systematically studied by Bowditch in higher-dimensional real hyperbolic spaces in \cite{B1}, and more generally in pinched Hadamard manifolds in \cite{B2}. In \cite{B1}, Bowditch labelled the five equivalent formulations of the definition of geometrical finiteness (GF1)-(GF5), with (GF3) corresponding to the original notion. He then showed in \cite{B2} that the four other formulations, now labelled F1,F2,F4, and F5, remain equivalent in the more general setting (but not the original one). The most convenient for our purposes will be condition F5, which we now recall. 

Let as above $X$ be a negatively curved symmetric space and $\Gamma$ a discrete subgroup of ${\rm Isom}(X)$. The \emph{convex hull} ${\rm Hull}(\Gamma)$ of $\Gamma$ in $X$ is the convex hull of the limit set $\Lambda(\Gamma)$, more precisely the smallest convex subset of $X$ whose closure in $\overline{X}=X \cup \partial_\infty X$ contains $\Lambda(\Gamma)$. This is invariant under the action of $\Gamma$, and the \emph{convex core} ${\rm Core}(\Gamma)$ of $\Gamma$ in $X$ is defined as the quotient of ${\rm Hull}(\Gamma)$ under the action of $\Gamma$. 

{\bf Definition:} We say that $\Gamma$ satisfies condition F5 if (a) for some $\varepsilon >0$, the tubular neighborhood $N_\varepsilon({\rm Core}(\Gamma))$ in $X/\Gamma$ has finite volume, and (b)  there is a bound on the orders of the finite subgroups of $\Gamma$.

\begin{prop}\label{geominf}
  The hybrid $H(3) < {\rm Isom}({\rm H}^2_\C)$ is geometrically infinite.
\end{prop}

\Pf We show that $H(3)$ does not satisfy condition F5. By Proposition~\ref{limitsets}, $\Lambda(H(3))=\partial_\infty {\rm H^2_\C}$, hence ${\rm Hull}(H(3))={\rm H^2_\C}$. Now by Theorem~\ref{d=3main}, $H(3)$ has infinite index in a lattice, therefore it acts on  ${\rm H^2_\C}$ with infinite covolume, in other words ${\rm Core}(H(3))$ has infinite volume hence so does any of its tubular neighborhoods. \EPf

\raggedright

\begin{flushleft}
  \textsc{Julien Paupert, Joseph Wells\\
   School of Mathematical and Statistical Sciences, Arizona State University}\\
       \verb|paupert@asu.edu, jswells@asu.edu|
\end{flushleft}

\frenchspacing
\begin{thebibliography}{ZZ99}
  
\bibitem[B1]{B1} B.H. Bowditch; {\sl Geometrical finiteness for hyperbolic groups}. J. Funct. Anal. {\bf 113} (1993), no. 2, 245--317. 
  
\bibitem[B2]{B2} B.H. Bowditch; {\sl Geometrical finiteness with variable negative curvature}. Duke Math. J. {\bf 77} (1995), no. 1, 229--274.
  
\bibitem[BM]{BM} Y. Benoist, S. Miquel; {\sl Arithmeticity of discrete subgroups
  containing horospherical lattices}. arXiv:1805.00045. 
  
\bibitem[BO]{BO} Y. Benoist, H. Oh; {\sl Discrete subgroups of ${\rm SL}_3(\R)$ generated by triangular matrices}. Int. Math. Res. Not. IMRN 2010, no. 4, 619--632.

\bibitem[CG]{CG} S.~Chen, L.~Greenberg; {\sl Hyperbolic spaces}, in Contributions to Analysis. Academic Press, New York (1974), 49--87.

\bibitem[DM]{DM} P.~Deligne, G.D.~Mostow; {\sl Monodromy of hypergeometric
functions and non-lattice integral monodromy}. Publ.~Math.~I.H.E.S.\ {\bf 63} (1986), 5--89.

\bibitem[D]{D} M. Deraux; {\sl A new non-arithmetic lattices in PU(3,1)} arXiv:1710.04463.

\bibitem[DF]{DF} M. Deraux, E. Falbel;  {\sl Complex hyperbolic geometry of the figure eight knot}. Geom. Topol. {\bf 19} (2015), 237--293. 

\bibitem[DPP1]{DPP1} M. Deraux, J.R.~Parker, J.~Paupert; {\sl New non-arithmetic complex hyperbolic lattices}. Invent. Math. {\bf 203} (2016), 681--771.

\bibitem[DPP2]{DPP2} M. Deraux, J.R.~Parker, J.~Paupert; {\sl On commensurability classes of non-arithmetic complex hyperbolic lattices}. arXiv:1611.00330.

\bibitem[F]{F} E. Falbel; {\sl A spherical CR structure on the complement of the figure eight knot with discrete holonomy}. J. Differential Geom. {\bf 79} (2008), no. 1, 69--110.

\bibitem[FFP]{FFP} E. Falbel, G. Francsics, J. R. Parker; {\sl The geometry of the Gauss-Picard modular group}. Math. Ann. {\bf 349} (2011), no. 2, 459--508. 

\bibitem[FP]{FP} E. Falbel, J. R. Parker; {\sl The geometry of the Eisenstein-Picard modular group}. Duke Math. J. {\bf 131} (2006), no. 2, 249--289.

\bibitem[FW]{FW} E. Falbel, J. Wang; {\sl Branched spherical CR structures on the complement of the figure eight knot}. Michigan Math. J. {\bf 63} (2014), no. 3, 635--667.

\bibitem[G]{G} W.M.~Goldman; Complex Hyperbolic Geometry.
Oxford Mathematical Monographs. Oxford University Press (1999).

\bibitem[GPS]{GPS} M.~Gromov, I.~Piatetski-Shapiro; {\sl Non-arithmetic groups in Lobachevsky spaces}. Publ. Math. IHES {\bf 66} (1987), 93--103. 

\bibitem[K1]{K1} M. Kapovich; {\sl On normal subgroups in the fundamental groups of complex surfaces}. arXiv: 
 9808085.
\bibitem[K2]{K2} M. Kapovich; {\sl Noncoherence of arithmetic hyperbolic lattices}. Geom. Topol. {\bf 17}  (2013), no. 1, 39--71.

\bibitem[K3]{K3} M. Kapovich; Hyperbolic manifolds and discrete groups. Progress in Mathematics {\bf 183}. Birkh\"{a}user (2001).

\bibitem[MP]{MP} A. Mark, J. Paupert; {\sl Presentations for cusped arithmetic hyperbolic lattices}. arXiv:1709.06691. 

\bibitem[M1]{M1} G.D.~Mostow; {\sl On a remarkable class of polyhedra
in complex hyperbolic space}. Pacific J.\ Maths.\ {\bf 86} (1980),
171--276.

\bibitem[M2]{M2} G.D.~Mostow; {\sl Generalized Picard lattices arising
from half-integral conditions}. Publ.~Math.~I.H.E.S.\ {\bf 63} (1986),
91--106.

\bibitem[O]{O} H. Oh; {\sl Discrete subgroups generated by lattices in opposite horospherical subgroups}. J. Algebra {\bf 203} (1998), no. 2, 621--676.
 
 \bibitem[Par1]{Par0} J.R.~Parker; {\sl Cone metrics on the sphere and Livn\'e's lattices}. Acta Math. {\bf 196}  (2006) 1--64.
 
\bibitem[Par2]{Par1} J.R.~Parker; {\sl Complex hyperbolic lattices}, in Discrete Groups and Geometric Structures. Contemp. Math. {\bf 501} AMS (2009), 1--42.
 
\bibitem[Par3]{Par2} J.R. Parker; Complex Hyperbolic Kleinian Groups. To appear.

\bibitem[Pau]{Pau} J. Paupert; {\sl Non-discrete hybrids in \SU(2,1)}. Geom. Dedicata {\bf 157} (2012), 259--268.

\bibitem[PW]{PW} J. Paupert, P. Will; {\sl Real reflections, commutators and cross-ratios in complex hyperbolic space}. Groups Geom. Dyn. {\bf 11} (2017), 311--352.

\bibitem[S]{S} P. Sarnak; {\sl Notes on thin matrix groups}, in Thin groups and superstrong approximation, 343--362,
Math. Sci. Res. Inst. Publ. {\bf 61}, Cambridge Univ. Press, Cambridge, 2014. 

\end{thebibliography}
\end{document}